\documentclass[11pt]{amsart}
\usepackage [T1]{fontenc}
\usepackage [latin1]{inputenc}
\usepackage{amssymb,amsmath,amsthm}
\usepackage{verbatim}
\usepackage{graphicx}
\usepackage{epsfig}

\usepackage{color, soul}
\usepackage[all]{xy}

\definecolor{dgreen}{cmyk}{1,.5,1,.2}

\newtheorem{theorem}{Theorem}[section]
\newtheorem{definition}[theorem]{Definition}

\newtheorem{example}[theorem]{Example}
\newtheorem{lemma}[theorem]{Lemma}
\newtheorem{proposition}[theorem]{Proposition}
\newtheorem{corollary}[theorem]{Corollary}
\newtheorem{remark}[theorem]{Remark}

\def\s{\backslash}
\begin{document}

\title[Unconditionality in tensor products]{Unconditionality in tensor products and ideals of polynomials, multilinear forms and operators}

%
%

\author{Daniel Carando}

\author{Daniel Galicer}

\address{Departamento de Matem\'{a}tica, Facultad de Ciencias Exactas y Naturales, Universidad de Buenos Aires,\\ Pab. I, Cdad Universitaria (1428)  Buenos Aires, Argentina and CONICET.} \email{dcarando@dm.uba.ar}
\email{dgalicer@dm.uba.ar}

\begin{abstract} We study tensor norms that destroy unconditionality in the following sense: for every Banach space $E$ with unconditional basis, the $n$-fold tensor product of $E$ (with the corresponding tensor norm) does not have unconditional basis. We establish an easy criterion to check weather a tensor norm destroys unconditionality or not.  Using this test we get that all injective and projective tensor norms different from $\varepsilon$ and $\pi$ destroy unconditionality, both in full and symmetric tensor products. We present applications to polynomial ideals: we show that many usual polynomial ideals never enjoy the Gordon-Lewis property. We also consider the unconditionality of the monomial basic sequence. Analogous problems for multilinear and operator ideals are addressed.
\end{abstract}

\maketitle

\section*{Introduction}
There has been a great interest on the study of unconditionality in tensor products of Banach spaces and, more recently, in spaces of polynomials and multilinear forms. As a rather uncomplete reference, we can mention \cite{DeDiGaMa01,DefKal05,DefPerGar08,PerVill04,PerVill05(unconditional),Pisier78,Schutt}.
A fundamental result obtained by Schütt~\cite{Schutt} and independently by Pisier~\cite{Pisier78} (with additional assumptions) simplified the study of unconditionality in tensor products: in order to know if a tensor product of Banach spaces with unconditional basis have also unconditional basis, just look at the monomials. The extension of these results to symmetric tensor norms (of any degree $n$) was probably motivated by the so called Dineen's problem or conjecture. In his book \cite{Din99}, Sean Dineen asked the following question: if the dual of a Banach space $E$ has an unconditional basis, can the space of $n$-homogeneous polynomials have unconditional basis? He conjectured a negative answer. Defant, Díaz, García and Maestre~\cite{DeDiGaMa01} developed the symmetric $n$-fold versions of Pisier and Schütt's work and, also,  obtained asymptotic estimates of the unconditionality constants of the monomial basis for spaces $\ell_p^m$. As a result, they made clear that a counterexample to Dineen's conjecture should be very hard to find. Finally, Defant and Kalton \cite{DefKal05} showed that if $E$ has unconditional basis, then the space of polynomials on $E$ cannot have unconditional basis. Defant and Kalton's result is based on a sort of dichotomy that they managed to establish: the space of polynomials either lacks the Gordon-Lewis property or is not separable. Therefore, should the space of polynomials have a basis, this cannot be unconditional.

On the other hand, in \cite{PerVill04} Pérez-García and Villanueva illustrated the bad behavior of many tensor norms with unconditionality. They showed, for example, than no natural tensor norm (in the sense of Grothendieck) preserve unconditionality: for any natural 2-fold tensor norm, there exists a Banach space with unconditional basis whose tensor product fails to have the Gordon-Lewis property.

Motivated by these results, we investigate when a tensor norm (of any degree, and either on the full or
on the symmetric tensor product) destroys unconditionality in the sense that, for \textbf{every} Banach space $E$ with unconditional basis, the corresponding tensor product has not unconditional basis.

We establish a simple criterion to check weather a tensor norm destroys unconditionality or not. With this we obtain that every injective and every projective tensor norm (other than $\varepsilon $ and $\pi$) destroys unconditionality. In particular, every non trivial symmetric natural norm (see the next section) destroys unconditionality.

In \cite{CarDim07,CaGa10} some differences between the $n=2$ and $n\ge 3$ cases were given. We present more evidence in which this cases are dissimilar: for $n=2$ the only natural tensor norms that destroy unconditionality are symmetric and for $n\ge 3$ there are non-symmetric natural tensor norms that destroy unconditionality. The contrasting situation between the $n=2$ and $n\ge 3$ cases is again exhibited for $n$-linear forms defined on the product of $n$ different spaces, as well as in tensor products of different spaces.

We also study unconditionality in ideals of polynomials and multilinear forms. We show that there are ideals $\mathcal{Q}^n$ of $n$-homogeneous polynomials such that, for every Banach space $E$ with unconditional basis, the space $\mathcal{Q}^n(E)$ lacks the Gordon-Lewis property. 
Among these ideals we have the $r$-integral, $r$-dominated, extendible and $r$-factorable polynomials. For the last three examples we even get
that the monomial basic sequence is never unconditional.

Note that the behavior of these ideals with unconditionality is more drastic than the ideal of all continuous polynomials $\mathcal{P}^n$, since $\mathcal{P}^n(E)$ can have the Gordon-Lewis property (for example, if $E=\ell_1$). We present another example of a maximal Banach polynomial ideal $\mathcal{Q}^n$ with the same property that Defant and Kalton showed for $\mathcal{P}^n$: $\mathcal{Q}^n(E)$ never has unconditional basis, but it may enjoy the Gordon-Lewis property. In this cases, $\mathcal{Q}^n(E)$ is not separable.

We consider ideals of multilinear forms and ideals of operators, where some results have their analogous.

\smallskip

We refer to \cite{DefFlo93} for the theory of tensor norms and operator ideals, and to \cite{Flo97,Flo01,Flo01(extension),Flo02(On-ideals),Flo02(Ultrastability)} for symmetric tensor products and polynomial ideals.

\section{Preliminaries}

For a natural number $n$, a full tensor norm $\alpha$ of order $n$ assigns to every $n$-tuple of Banach spaces $(E_1, \dots, E_n)$ a norm $\alpha \big(\; . \; ; \otimes_{i=1}^n E_i \big)$ on the $n$-fold (full) tensor product $\otimes_{i=1}^n E_i$ such that

\begin{enumerate}
\item $\varepsilon \leq \alpha \leq \pi$ on $\otimes_{i=1}^n E_i$.
\item $\| \otimes_{i=1}^n T_i :  \big( \otimes_{i=1}^n E_i, \alpha \big) \to \big( \otimes_{i=1}^n F_i, \alpha \big) \|  \leq \|T_1\| \dots \|T_n\|$ for each set of operator $T_i~\in~\mathcal{L}(E_i, F_i)$, $i=1, \dots, n$.
\end{enumerate}

Here, $\varepsilon$ and $\pi$ denote the injective and projective tensor norms respectively.

We say that $\alpha$ is finitely generated if for all Banach spaces $E_i $  and each $z$ in $\otimes_{i=1}^n E_i$ we have
$$ \alpha ( z, \otimes_{i=1}^n E_i ) : = \inf \{ \alpha (z, \otimes_{i=1}^n M_n ) : z \in \otimes_{i=1}^n M_i \},$$
the infimum being taken over all $n$-tuples $M_1, \dots, M_n$ of finite dimensional subspaces of $E_1, \dots, E_n$ respectively whose tensor product contains $z$.

We often call these tensor norms ``full tensor norms'', in the sense that they are defined on the full tensor product, to distinguish them from the s-tensor norms, that are defined on symmetric tensor products.

We say that  $\beta$ is an s-tensor norm  of order $n$ if $\beta$ assigns to each Banach space $E$ a norm $\beta \big(\; . \;; \otimes^{n,s} E \big)$ on the $n$-fold symmetric tensor product $\otimes^{n,s} E$ such that
\begin{enumerate}
\item $\varepsilon_s \leq \beta \leq \pi_S$ on $\otimes^{n,s} E$.
\item $\| \otimes^{n,s} T :  \big( \otimes^{n,s} E, \beta \big) \to \big( \otimes^{n,s} F, \beta \big) \| \leq \|T\|^n$ for each operator $T~\in~\mathcal{L}(E, F)$.
\end{enumerate}

$\beta$ is called finitely generated if for all $E \in  BAN$ and $z \in \otimes^{n,s} E$
$$ \beta (z, \otimes^{n,s}E) = \inf \{ \alpha(z, \otimes^{n,s}M) : M \in FIN(E), z \in \otimes^{n,s}M \}.$$

In both cases condition $(2)$ will be referred to as the ``metric mapping property''. Also, all the full tensor norms (or s-tensor norms) that are not equivalent to $\varepsilon$ (or $\varepsilon_s$) nor $\pi$ (or $\pi_s$) will be referred to as \textbf{nontrivial}.

Throughout the article, we will assume that all tensor norms are finitely generated.

If $\alpha$ is a full tensor norm of order $n$, then the dual tensor norm $\alpha'$ is defined on FIN (the class of finite dimensional Banach spaces) by

$$  \big( \otimes_{i=1}^n M_i, \alpha' \big) :\overset 1 = [\big( \otimes_{i=1}^n M_i', \alpha \big)]'$$
and on BAN (the class of all Banach spaces) by
$$ \alpha' ( z, \otimes_{i=1}^n E_i ) : = \inf \{ \alpha' (z, \otimes_{i=1}^n M_n ) : z \in \otimes_{i=1}^n M_i \},$$
the infimum being taken over all $n$-tuples $M_1, \dots, M_n$ of finite dimensional subspaces of $E_1, \dots, E_n$ respectively whose tensor product contains $z$.

Analogously, for $\beta$ an s-tensor norm of order $n$, its dual tensor norm $\beta'$ is defined on FIN by

$$  \big( \otimes^{n,s} M, \beta' \big) :\overset 1 = [\big( \otimes^{n,s} M', \beta \big)]'$$
and extended to BAN as before.

For $\alpha$, a full tensor norm of order $n$, we will denote $\underline{\alpha}$ the full tensor norm of order $n-1$ given by $$\underline{\alpha}(z, \otimes_{i=1}^{n-1} E_i) := \alpha(z \otimes 1,  E_1 \otimes \dots \otimes  E_{n-1} \otimes \mathbb{C}),$$ where  $z \otimes 1 := \sum_{i=1}^m x_1^i \otimes \dots x_n^i \otimes 1$, for $z = \sum_{i=1}^m x_1^i \otimes \dots x_n^i$ (this definition can be seen as dual to some ideas on  \cite{BotBraJunPel06} and \cite{CarDimMurN}).

\medskip

Using the metric mapping property and the definition of the operation $\underline{(\cdot)}$ we get the following remark, which will be used throughout the article:

\begin{remark}\label{bajar_de_n_a_2}\rm
Let $E_1, \dots, E_n$ Banach spaces, $x_j \in B_{E_j}$ ($j=3,\dots,n$) and $\alpha$ a full tensor norm of order $n$. Then $\big( E_1 \otimes E_2 \otimes [x_3] \otimes \dots \otimes [x_n], \alpha \big)$ is a complemented subspace of $\big( E_1  \dots \otimes E_n, \alpha \big)$ and this space is isometrically isomorphic to $\big(E_1 \otimes E_2, \widetilde{\alpha} \big)$, where $\widetilde{\alpha}$ is the 2-fold tensor norm which comes from applying  $n-2$ times the operation $\underline{(\cdot)}$ to the norm $\alpha$.
\end{remark}

Let $\alpha$ be a full tensor norm of order $n$.
We will say that $\alpha$ is projective if, whenever $P_i: E_i \to F_i$  are quotient maps ($i=1 \dots n$), the tensor product operator
$$ \otimes_{i=1}^n P_i : \big( \otimes_{i=1}^n E_i, \alpha \big) \to \big( \otimes_{i=1}^n F_i, \alpha \big),$$
is also a quotient map.

On the other hand, we will say that $\alpha$ is injective if, whenever $I_i: E_i \to F_i$  are isometric embeddings ($i=1 \dots n$), the tensor product operator
$$ \otimes_{i=1}^n I_i : \big( \otimes_{i=1}^n E_i, \alpha \big) \to \big( \otimes_{i=1}^n F_i, \alpha \big),$$
is an isometric embedding.

The projective and injective associates (or hulls) of $\alpha$ will be  denoted, by extrapolation of the 2-fold case, as $\s \alpha /$ and $/ \alpha \s$ respectively. The projective associate of $\alpha$ will be the (unique) smallest projective tensor norm greater than $\alpha$. Following \cite[Theorem 20.6.]{DefFlo93} we have:
$$ \big( \otimes_{i=1}^n \ell_1(B_{E_i}),  \alpha  \big) \overset 1 \twoheadrightarrow \big( \otimes_{i=1}^n  E_i,   \s \alpha /  \big).$$
The injective associate of $\alpha$ will be the (unique) greatest injective tensor norm smaller than $\alpha$.  As in \cite[Theorem 20.7.]{DefFlo93} we get,
$$ \big( \otimes_{i=1}^n E_i , / \alpha \s \big) \overset 1 \hookrightarrow  \big( \otimes_{i=1}^n \ell_{\infty}(B_{E_i'}),  \alpha  \big).$$
Note that in our notation, the symbols ``$\s$'' and ``$/$'' by themselves lose their original meanings, as well as the left and right sides of $\alpha$.

With this, an $n$-linear form $A$ belongs to $\big( \otimes_{i=1}^n  E_i,   \s \alpha /  \big)'$ if and only if $A \circ ( P_{E_1}, \dots, P_{E_n}) \in \big( \otimes_{i=1}^n  \ell_1(B_{E_i}),   \alpha   \big)'$ where $P_{E_1} : \ell_1(B_E) \twoheadrightarrow E$ stands for the canonical quotient map.
Moreover, $$\|A\|_{\big( \otimes_{i=1}^n  E_i,   \s \alpha /  \big)'}= \|A \circ ( P_{E_1}, \dots, P_{E_n}) \|_{\big( \otimes_{i=1}^n  E_i,   \alpha  \big)'}.$$
On the other hand, an $n$-linear form $A$ is  $\big( \otimes_{i=1}^n  E_i,   / \alpha \s  \big)'$ if it has an extension to $\ell_{\infty}(B_{E_1'})\times\dots\times \ell_{\infty}(B_{E_n'})$ that is $ \big( \otimes_{i=1}^n  \ell_{\infty}(B_{E_i'}),   \alpha   \big)'$. Moreover, the norm of $A$ in $\big( \otimes_{i=1}^n  E_i,   / \alpha \s  \big)'$ is the infimum of the norms in $\big( \otimes_{i=1}^n  \ell_{\infty}(B_{E_i'}),    \alpha   \big)'$ of all such extensions.

The projective and injective associates for an s-tensor norm $\beta$ can be defined in a similar way:
$$ \big( \otimes^{n,s} \ell_1(B_{E}),  \beta  \big) \overset 1 \twoheadrightarrow \big( \otimes^{n,s}  E,   \s \beta /  \big).$$
$$ \big( \otimes_{i=1}^n E, / \beta \s \big) \overset 1 \hookrightarrow  \big( \otimes^{n,s} \ell_{\infty}(B_{E'}),  \beta \big).$$

The description of the $n$-homogeneous polynomial $Q$ belonging to $\big( \otimes^{n,s}  E,   \s \beta /  \big)'$ or to $\big( \otimes^{n,s}  E,   / \beta \s  \big)'$ is analogous to that for multilinear forms.

Following the ideas of \cite[Proposition 20.10.]{DefFlo93}, for a full tensor norm $\alpha$ and an s-tensor norm $\beta$, we have the following duality relations
$$ (/ \alpha \s)' = \s \alpha ' /, \; \; \; (\s \alpha /)' = /\alpha '\s, \; \; \;   (/ \beta \s)' = \s \beta' /, \; \; \;   (\s \beta /)' = / \beta' \s.$$

Note that one could have defined a first, second, up to $n$th injective and projective associates (and any combination of them), in the spirit of the right and left injective and projective associates for 2-fold tensor products \cite[Section 20]{DefFlo93}. However, the notation would be rather uncomfortable and we will only use these associates in just one example. Therefore, we will only introduce the notation for that particular case when necessary.

\bigskip
In  \cite{CaGa10} we have introduced and studied natural full symmetric tensor norm (of any order) in the spirit of the work of Grothendieck. Those are obtained from $\pi$ with a finite number of the operations $\setminus  \ /$, $/ \ \setminus$, $'$. There are exactly four natural full symmetric tensor norms of order 2 \cite[Section 27]{DefFlo93}, namely, $\pi$, $\varepsilon$, $/\pi\s \sim w_2$ and $\s \varepsilon / \sim w_2'$.
But for $n\ge 3$ we actually have 6 different ones. They can be arranged in the following way:
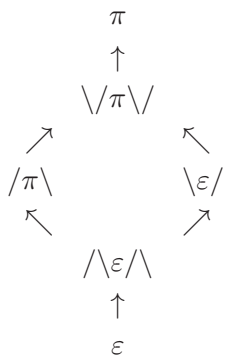
\begin{equation*}
\begin{array}{rcl}
 & \pi &  \\
 &  \uparrow &  \\
 & \s / \pi \s / &  \\
 \nearrow & &\nwarrow  \\
/ \pi \s &  & \s \varepsilon / \\
\nwarrow & &  \nearrow \\
& / \s \varepsilon / \s & \\
 &  \uparrow &  \\
& \varepsilon & \\
\end{array}
\end{equation*}
where $\alpha \to \gamma$ means that $\gamma$ dominates $\alpha$. And there are no other dominations than those showed in the scheme. We will discuss unconditionality for natural norms in the next section.

\section{Destruction of unconditionality}

In their fundamental paper \cite{GordonLewis1974} from 1974, Gordon and Lewis showed that spaces of operator between infinite dimensional Banach spaces lack a ``reasonable unconditional structure'', particularly the space $\mathcal{L}(\ell_2)$ of all operators on the Hilbert space $\ell_2$. Their key was to prove the so-called Gordon-Lewis inequality which estimates the unconditional basis constant by its Gordon Lewis constant.

Recall that a Banach space $E$ has the Gordon-Lewis property if every absolutely summing operator $T: E \to \ell_2$ is $1$-factorable (i.e. allows a factorization $T: E \overset R \to L_1(\mu) \overset S \to \ell_2$). In this case, there is a constant $C \geq 0$ such that for all $T: E \to \ell_2$,
$$\gamma_1(T) : = \inf \|R\| \|S\| \leq C \pi_1(T),$$ and the best constant $C$ is called the Gordon-Lewis constant of $E$ and denoted by $gl(E)$.

These ideas were taken by Pisier \cite{Pisier78} and Schütt \cite{Schutt} to give a deep study of unconditionality in tensor product of Banach spaces. They showed (independently) that for any full tensor norm $\alpha$ on the tensor product $E \otimes F$ of two Banach spaces with unconditional basis $(e_i)$ and $(f_j)$, respectively, the monomials $(e_i \otimes f_j)_{i,j}$ form an unconditional basis if and only if $E \widetilde{\otimes}_{\alpha} F$ has unconditional basis if and only if $E \widetilde{\otimes}_{\alpha} F$ has the Gordon Lewis property. This was generalized by Defant, Díaz, García y Maestre in \cite{DeDiGaMa01} to the $n$-fold case.

\begin{theorem} \cite[Remark 1]{DeDiGaMa01}\label{full-Gl_incond}
Let $E_1, \dots, E_n$ be a finite sequence of Banach spaces with 1-unconditional basis $(e_k^j)$. Then for each full tensor norm $\alpha$,
$$\chi_{mon} \Big( \big( \widetilde{\otimes}_{j=1}^n E_j, \alpha \big) \Big) \leq 2^{n+1} gl \Big( \big( \widetilde{\otimes}_{j=1}^n E_j, \alpha \big) \Big),$$
where $\chi_{mon}$ stands for the unconditional basis constant of to the monomial basis.
\end{theorem}
 In particular this shows that $\big( \widetilde{\otimes}_{j=1}^n E_j, \alpha \big)$ has unconditional basis if and only if the monomials form an unconditional basis of $\big( \widetilde{\otimes}_{j=1}^n E_j, \alpha \big)$ if and only if $\big( \widetilde{\otimes}_{j=1}^n E_j, \alpha \big)$ has the Gordon-Lewis property, whenever $E_j$ has unconditional basis.

A similar result also holds for the symmetric $n$-fold tensor product of a Banach space with unconditional basis:

\begin{theorem} \cite[Corollary 1.]{DeDiGaMa01}\label{symm-Gl-Incond}
Let $E$ be a Banach space with unconditional basis $(e_j)$. Then for each s-tensor norm $\beta$ of order $n$, the following are equivalent:
\begin{enumerate}
\item The monomials of order $n$ with respect to $(e_j)$ form an unconditional basis of $\big( \widetilde{\otimes}^{n,s} E , \beta \big)$
\item $\big( \widetilde{\otimes}^{n,s} E , \beta \big)$ has unconditional basis
\item $\big( \widetilde{\otimes}^{n,s} E , \beta \big)$ has the Gordon-Lewis property.
\end{enumerate}
\end{theorem}

An interesting result due to Pérez-García and Villanueva \cite[Proposition 2.3]{PerVill04} is that, if $\big( \widetilde{\otimes}_{i=1}^n c_0, \alpha \big)$ has unconditional basis, then $\alpha$ has to coincide (up to constants) with the injective norm $\varepsilon$ on $\otimes_{i=1}^n c_0$. On the other hand, if the tensor product $\big( \widetilde{\otimes}_{i=1}^n \ell_1, \alpha \big)$ has unconditional basis then $\alpha$ has to be equivalent to the projective norm $\pi$ on $\otimes_{i=1}^n \ell_1$ \cite[Proposition 2.6]{PerVill04}.

A similar statement holds when considering Hilbert spaces \cite[Theorem 2.5.]{PerVill05(unconditional)}. More precisely, if $\big( \widetilde{\otimes}_{i=1}^n \ell_2, \alpha \big)$ has unconditional basis then $\alpha$ has to coincide with the Hilbert-Schmidt norm $\sigma_2$ (again, up to constants).
It is important to remark that if $\big( \widetilde{\otimes}_{i=1}^n c_0, \alpha \big)$ and $\big( \widetilde{\otimes}_{i=1}^n \ell_1, \alpha \big)$ have both unconditional basis this imply that on the tensor product of Hilbert spaces $\alpha$ equals the Hilbert-Schmidt norm (we get this from \cite[Proposition 2.7.]{PerVill04} and \cite[Theorem 4.2.]{Per04(inclusion)})

A consequence of this results is the following, also due to Pérez-García and Villanueva \cite[Theorem 1.1]{PerVill04}, is that \textit{no natural 2-fold tensor norm preserves unconditionality}: if $\alpha$ is a 2-fold natural tensor norm, there exists a Banach space $E$ with unconditional basis such that $E \widetilde{\otimes}_{\alpha} E$ has not unconditional basis.

It is not hard to see that the same holds for natural tensor norms of higher order. Note that for a tensor norm, ``not preserving'' unconditionality means that there is some space with unconditional basis such that its tensor product lacks of it. We will see that many of the natural tensor norms have a more drastic behavior: they destroy unconditionality in the following sense:

\begin{definition} We will say that a full tensor norm $\alpha$ \textbf{destroys unconditionality} if the tensor product $\big( \widetilde{\otimes}^n E, {\alpha} \big)$ does not have unconditional basis for any Banach space $E$ with unconditional basis.
\end{definition}
As mentioned before, in order that a full tensor norm $\alpha$ preserve unconditionality it is necessary for $\alpha$ to be equivalent to $\varepsilon$, $\sigma_2$ and $\pi$ in $\otimes^n c_0$, $\otimes^n \ell_2$, $\otimes^n \ell_1$ respectively. If none of these conditions are satisfied, we have just the opposite: $\alpha$ destroys unconditionality.

\bigskip

\begin{theorem} \textbf{Destruction Test:} \label{Destruction test}
A full tensor norm $\alpha$ destroys unconditionality if and only if $\alpha$ is not equivalent to $\varepsilon$, $\sigma_2$ and $\pi$ on $\otimes^n c_0$, $\otimes^n \ell_2$ and $\otimes^n \ell_1$ respectively.
\end{theorem}

To prove this we will need a result of Tzafriri:

\begin{theorem}\cite{Tzafriri} \label{tzafriri}
Let $E$ be a Banach space with unconditional basis then $E$ contains uniformly complemented at least one of the three sequences $(\ell_p^m)_{m=1}^\infty$ with $p \in \{1,2, \infty\}$.
\end{theorem}

\begin{proof} (of Theorem~\ref{Destruction test}) It is clear that a tensor norm that destroys unconditionality cannot enjoy any of the three equivalences in the statement. Conversely, suppose that $\alpha$ is not equivalent to $\varepsilon$, $\sigma_2$ and $\pi$ on $\otimes^n c_0$, $\otimes^n \ell_2$ and $\otimes^n \ell_1$ respectively. Let us see that if $E$ be a Banach space with unconditional basis, then $\big( \widetilde{\otimes}^n E, \alpha \big)$ cannot not have the Gordon-Lewis property.
By Theorem~\ref{tzafriri} we know that $E$ contains an uniformly complemented sequence of $(\ell_p^m)_{m=1}^\infty$ for $p=1,2$ or $\infty$. So, fixed such $p$, there is a constant $K>0$ such that
$$
gl \Big( \big(\otimes^n \ell_p^m , \alpha \big) \Big) \leq K gl \Big( \big(\widetilde{\otimes}^n E , \alpha \big) \Big),
$$
for every $m$.
If $gl \Big( \big(\widetilde{\otimes}^n E , \alpha \big) \Big)$ is finite then, by Theorem~\ref{full-Gl_incond},  $$\chi_{mon}\Big( \big( \widetilde{\otimes}^n \ell_p , \alpha \big) \Big) = \sup_{m} \chi_{mon}\Big( \big( \otimes^n \ell_p^m , \alpha \big) \Big) < \infty, \mbox{if $p=1$ or $2$},$$  or $$\chi_{mon}\Big( \big( \widetilde{\otimes}^n c_0 , \alpha \big) \Big)= \sup_{m} \chi_{mon}\Big( \big( \otimes^n \ell_\infty^m , \alpha \big) \Big) < \infty \mbox{ if $p=\infty$}.$$
This implies that either $\big( \widetilde{\otimes}^n \ell_1 , \alpha \big)$ or $\big( \widetilde{\otimes}^n \ell_2 , \alpha \big)$  or  $\big( \widetilde{\otimes}^n c_0 , \alpha \big)$ has unconditional basis. Now using Pérez-García and Villanueva results \cite[Propositions 2.3 and 2.6]{PerVill04} and \cite[Theorem 2.5]{PerVill05(unconditional)} we have that,  either  $\alpha \sim \pi$ on $\otimes^n \ell_1$ or $\alpha \sim \sigma_2$ on $\otimes^n \ell_2$ or  $\alpha \sim \varepsilon$ on $\otimes^n c_0$, which leads us to a contradiction. Therefore, $gl \Big( \big(\widetilde{\otimes}^n E , \alpha \big) \Big)$ is infinite and the statement is proved. \end{proof}

As a simple consequence of the test we have that, a full tensor norm $\alpha$ destroys unconditionality if and only if $\alpha'$ destroys unconditionality.



\medskip
We will show that injective or projective tensor norms other than $\varepsilon$ and $\pi$ destroy unconditionality.

First, note that from \cite[Proposition 3.1]{CarDimSev06} (and its proof), we can see that if $S$ is a diagonal extendible multilinear form on $\ell_p$ ($2\le p\le\infty$), then $S$ is nuclear and
\begin{equation}\label{extendible-nuclear} \|S\|_{\mathcal{N}} \leq C \|S \|_e.\end{equation} The definition of nuclear and extendible multilinear forms can also found in \cite{CarDimSev06}  (in the Section~\ref{section-ideals} we present the definition of extendible polynomials, which is analogous). Just for completeness, extendible multilinear forms are exactly those that are $/ \pi \s$-continuous, and nuclear multilinear forms are $\varepsilon$-continuous.

If $T$ is any multilinear form on $\ell_p$, we denote by $D(T)$ the multilinear form obtained from $T$ setting to zero all the coefficients outside the diagonal (see \cite{CarDimSev07} for details).

\begin{lemma} \label{eta-epsilon}
Let $2 \leq p \leq \infty$. There exist a constant $K$ such that for every sequence of scalars $a_1, \dots, a_m$,
$$ / \pi \s ( \sum_{k=1}^m a_k e_k \otimes \dots \otimes e_k, \otimes^n \ell_p^m) \leq K \varepsilon(\sum_{k=1}^m a_k e_k \otimes \dots \otimes e_k, \otimes^n \ell_p^m)$$
\end{lemma}

\begin{proof}
Notice that
\begin{align*}
/ \pi \s \big( \sum_{k=1}^m a_k e_k \otimes \dots \otimes e_k, \otimes^n \ell_p^m \big)
& = \sup_{ \|T\|_e \leq 1} \big | \sum_{k=1}^m a_k T(e_k, \dots, e_k) \big| \\
& = \sup_{ \|T\|_e \leq 1} \big | \sum_{k=1}^m a_k D(T)(e_k, \dots, e_k) \big|, \\
& \leq \sup {\left\{\big| \sum_{k=1}^m a_k S (e_k, \dots, e_k) \big| : S \in \mathcal{L}^n( \ell_p^m) \; \mbox{diagonal} \; : \; \|S\|_e \leq 1 \right\}} ,
\end{align*}
where the last inequality is a consequence of the inequality $\|D(T)\|_e \leq \|T\|_e$ (\cite[Proposition~5.1.]{CarDimSev07}).
Now, using (\ref{extendible-nuclear}), we have
\begin{align*}
/ \pi \s \big( \sum_{k=1}^m a_k e_k \otimes \dots \otimes e_k, \otimes^n \ell_p^m \big) & \leq C^{-1}
\sup{\left \{ \big| \sum_{k=1}^m a_k S (e_k, \dots, e_k) \big|: S \in \mathcal{L}^n( \ell_p^m) \; \mbox{diagonal} \; : \; \|S\|_\mathcal{N} \leq 1 \right \}}  \\
& \le C^{-1} \varepsilon \big( \sum_{k=1}^m a_k e_k \otimes \dots \otimes e_k, \otimes^n \ell_p^m \big).
\end{align*}
\end{proof}

Now, what we are ready to show:

\begin{theorem}\label{inyectivas y proyectivas destruyen} Nontrivial injective and nontrivial projective tensor norms destroy unconditionality.
\end{theorem}

\begin{proof}
Let us see first that $/ \pi \s$ destroys unconditionality. By the Destruction Test (Theorem~\ref{Destruction test}) we need to show that  $/ \pi \s$ is not equivalent to $\varepsilon$, $\sigma_2$ and $\pi$ on $\otimes^n c_0$, $\otimes^n \ell_2$ and $\otimes^n \ell_1$ respectively.

The tensor norm $/ \pi \s$ is not equivalent to $\varepsilon$ on $\otimes^n c_0$: since $/ \pi \s = \pi$ on $\otimes^n c_0$, this would imply $\pi \sim \varepsilon$, which clearly false.

\smallskip

The tensor norm $/ \pi \s$ is not equivalent to $\sigma_2$ on $\otimes^n \ell_2$:
Lemma~\ref{eta-epsilon} states the existence of a constant $K$ such that:
$$ / \pi \s \big( \sum_{k=1}^m e_k \otimes \dots \otimes e_k, \otimes^n \ell_2 \big) \leq K \varepsilon\big( \sum_{k=1}^m e_k \otimes \dots \otimes e_k, \otimes^n \ell_2 \big) \leq K. $$
On the other hand, $$\sigma_2 \big(\sum_{k=1}^m e_k \otimes \dots \otimes e_k, \otimes^n \ell_2  \big) = m^{1/2},$$
so we are done. Note that this shows that
\begin{equation} \label{eq eta sigma 2}
\| id: \big( \otimes^n \ell_2^m, {/ \pi\s} \big) \longrightarrow \big( \otimes^n \ell_2^m, {\sigma_2} \big)\| \to \infty,
\end{equation}
as $m \to \infty$, a fact that will be used below.

The tensor norm $/ \pi \s$ is not equivalent to $\pi$ on $\otimes^n \ell_1$:
if it were,  every $n$-linear form on $\ell_1$ would be extendible, but this cannot happen (see, for example,  \cite[Corollary 12]{Carando01} ). Since $/ \pi \s \leq \pi$, this shows that
\begin{equation} \label{eq eta pi}
\| id: \big( \otimes^n \ell_1^m, {/ \pi\s} \big)\longrightarrow \big( \otimes^n \ell_1^m, {\pi} \big)\| \to \infty,
\end{equation}
as $m \to \infty$.

Thus, we have shown that $/ \pi \s$ destroys unconditionality. From Equations~(\ref{eq eta sigma 2}) and (\ref{eq eta pi}), if $\alpha$ is a tensor norm that is dominated by $/ \pi \s$, then it cannot be equivalent to $\pi$ or $\sigma_2$ on $\otimes^m \ell_1$ or $\otimes^m \ell_2$ respectively. If it is equivalent to $\varepsilon$ on $c_0$, we would have that $/ \alpha \s$ must be equivalent to $\varepsilon$ (on BAN). Therefore, the only (up to equivalences) injective tensor norm that does not destroy unconditionality is $\varepsilon$. By duality, a projective tensor norm that is not equivalent to $\pi$ must destroy unconditionality.

%
%
\end{proof}
%
%
%

The previous result asserts that nontrivial natural full-symmetric tensor norms destroy unconditionality. A natural question arises: what about the other (non-symmetric) natural norms? We know that none of them preserve unconditionality, but which of them destroy it? Again, the answer will depend on $n$ being 2 or greater:

\begin{remark}
\rm For $n=2$, $/ \pi \s$ and $\s \varepsilon /$ are the only natural norms that destroy unconditionality.
\end{remark}
\begin{proof}
We know that $/ \pi \s$ and $\s \varepsilon /$ destroy unconditionality and that  $\pi$ and $\varepsilon$ do not.

On the other hand, since $(/ \pi \s)/ \sim d_2$ is equivalent to  $\sigma_2$ in $\otimes^2 \ell_2$, we have that $(/ \pi \s)/ $ does not destroy unconditionality and, by duality, neither does  $\s(/ \pi \s) \sim g_2$.

By \cite[Corollary 3.2]{I.Schutt} we know that  $\Pi_1(\ell_2, \ell_2)$ has the Gordon-Lewis property. Therefore, $\varepsilon / = d_\infty$ cannot destroy unconditionality. Transposing and/or dualizing, neither do $\s \varepsilon = g_\infty$, $ \pi \s = d_\infty '$ or $ / \pi = g_\infty ' $.

If we show that $\s (/ \pi)= \s g_\infty '$ does not destroy unconditionality, we obtain the same conclusion for $(\pi \s )/ = d_\infty ' /$, $(\varepsilon / ) \s = d_\infty \s$ and $/ ( \s \varepsilon) = / g_\infty$ (again by duality and trasnposition)
Now, since $ \ell_\infty$ is injective, every operator from $\ell_1$ to  $ \ell_\infty$ is extendible. Therefore, $/ \pi$ and $\pi$ are equivalent on $\otimes^2 \ell_1$, which implies also the equivalence of  $\s (/ \pi)$ and $\pi$ on $\otimes^2 \ell_1$, and thus $\s (/ \pi)\sim \s g_\infty '$ does not destroy unconditionality, which ends the proof.
\end{proof}

We have just shown that, for $n=2$, nontrivial symmetric tensor norms are exactly those that destroy unconditionality. Let us see that for $n\ge 3$, there are non-symmetric natural tensor norms that destroy unconditionality. We have never defined nor introduced the the notation for non-symmetric natural tensor norms, but for the following examples, it is enough to say that $inj_k$ means to take injective associate in the $k$th place (e.g., for $n=2$, $inj_1\alpha$ is the left injective associate $/\alpha$).

\begin{example}
There are non-symmetric natural norms that destroy unconditionality.
\end{example}
Consider $\alpha =  inj_2 \; inj_1 \pi_n $. Note that $E \otimes_{/\pi\setminus} E$ is isometric to a complemented subspace of $\big( \otimes^n E, \alpha \big)$ for any Banach space $E$. Since $/\pi\setminus$ destroy unconditionality, it destroys the Gordon-Lewis property, and therefore so does $\alpha$.

\smallskip
It is not true that every natural tensor norm different from $\pi$ and $\varepsilon$ destroys unconditionality. For example, if we take $\alpha = inj_1 \pi$ we have
$$ (\ell_1 \otimes \ell_1 \otimes \ell_1, \alpha) \simeq (\ell_1 \otimes_{ / \pi_2} \ell_1) \otimes_{\pi_2} \ell_1 \simeq \big( \otimes^3 \ell_1,  \pi_3 \big).$$
Therefore,  $\alpha = inj_1 \pi$ does not destroy unconditionality.

\bigskip
Our original motivation was the unconditionality problem for spaces of polynomials (Dineen's problem), and so it was reasonable to consider tensor products of a single space. However, the question about unconditionality is interesting also in tensor products of different spaces. Moreover, we will see that in this case, there is a new difference between $n=2$ and $n\ge 3$. First we have this lemma (we present a proof at the end of this section):


\begin{lemma}\label{GL-chevet}
Let $\alpha$ be a 2-fold full injective norm. There exist a constant $C \geq 0$ such that $ m^{1/2} \leq C gl \big( \ell_1^m \otimes_{\alpha} \ell_2^m \big) \leq C m^{1/2} $ for every $m \in \mathbb{N}$. In particular, $gl \big( \ell_1^m \otimes_{\alpha} \ell_2^m \big)  \to \infty,$ as $m \to \infty$.
\end{lemma}

Now we can prove:

\begin{proposition} \label{Para mas de 3 arruina}
Fix $n\ge 3$ and let $\alpha$ be an $n$-fold natural full symmetric tensor norm other than $\pi$ or $\varepsilon$. If $E_1, \dots, E_n$ have unconditional bases, then  $\otimes_{\alpha}(E_1, \dots, E_n)$ does not have the Gordon-Lewis property (nor unconditional basis).
\end{proposition}

\begin{proof}
By the previous preposition we have that \begin{equation}\label{gl2tensor}gl \big( \ell_1^m \otimes_{/ \pi \s} \ell_2^n \big) \simeq gl \big( \ell_1^m \otimes_{/ \s \varepsilon / \s } \ell_2^n \big) = gl \big( \ell_{\infty}^m \otimes_{\s / \pi \s /} \ell_2^n \big) \simeq gl \big( \ell_{\infty}^m \otimes_{\s  \varepsilon /} \ell_2^n \big),\end{equation} and all go to infinity as $m$ goes to infinity (the constants involved in the equivalences do not depend on $m$). Let $\tilde \alpha$ as in Remark~\ref{bajar_de_n_a_2}, it is easy to show that $\tilde \alpha$ is the 2-fold natural analogous to $\alpha$, thus must be one of the tensor norms that appear in (\ref{gl2tensor}). Recall that nontrivial natural symmetric tensor norms destroy unconditionality, therefore $gl \big( \ell_p^m \otimes_{\tilde \alpha} \ell_p^m \big) \to \infty$ for $p \in \{ 1, 2, \infty \}$.

By Theorem~\ref{tzafriri}, since we have at least 3 spaces, two of them must contain, respectively, $\ell_p^n$'s and  $\ell_q^n$'s uniformly complemented, for $p$ and $q$ such that $gl \big(\ell_p^n \otimes_{\tilde \alpha} \ell_q^n \big) \to \infty$. Say $E_1$ and $E_2$ are those spaces. Observe that  $\ell_p^n \otimes_{\tilde \alpha} \ell_q^n = \big(\ell_p^n \otimes \ell_q^n \otimes [e_3] \otimes \dots [e_n], \alpha \big)$ are uniformly complemented in $\big( \widetilde{\otimes}^n E_i, \alpha \big)$ by Remark~\ref{bajar_de_n_a_2} and the proof is complete.
\end{proof}

With a similar proof the same result holds for $\alpha$ an $n$-fold injective (projective) full tensor norm such that $\tilde \alpha \not \sim \varepsilon$ ($\tilde \alpha \not \sim \pi$). It is important to note that the previous proposition is false for $n=2$. Indeed,
$ c_0~\otimes_{ / \pi \s}~\ell_2 =c_0 \otimes_{ \pi \s} \ell_2 = c_0 \otimes_{ d_\infty'} \ell_2$, so if we show that there exists $C>0$ such that
 $gl(\ell_\infty^m \otimes_{d_\infty'} \ell_2^m) \leq C$ for every $m$, we are done.
We have $$gl \big( \ell_\infty^m \otimes_{d_\infty'} \ell_2^m \big) = gl \big( \ell_1^m \otimes_{d_\infty} \ell_2^m \big) = gl \big( ( \ell_1^m \otimes_{d_\infty} \ell_2^m)' \big)= gl\big( \Pi_1(\ell_1^m,\ell_2^m) \big).$$
In \cite{I.Schutt}, I. Schütt showed that the last expression is uniformly bounded. This fact can be deduced easily in a different way. Indeed, by Grothendieck's Theorem (one of them!) we have that $\Pi_1 (\ell_1, \ell_2) = \mathcal{L}(\ell_1,\ell_2)$, then $gl\big( \Pi_1(\ell_1^m,\ell_2^m) \big) \asymp gl \big( \mathcal L(\ell_1^m,\ell_2^m) \big) = gl \big( \ell_{\infty}^m \otimes_{\varepsilon} \ell_2^m \big)$ where the equivalence constants are independent of $m$. Since $\chi\big( \ell_{\infty}^m \otimes_{\varepsilon} X \big)=1$ for every space such $\chi(X)=1$ we are done.

\bigskip

In \cite[Theorem 2.6]{PerVill05(unconditional)} it is proved that, if $\beta$ is a s-tensor norm such that $\big( \widetilde{\otimes}^{n,s} \ell_2, \beta \big)$ has an unconditional basis, then $\beta$ has to be equivalent to the Hilbert-Schmidt s-tensor norm.
The analogous result for $c_0$ and $\ell_1$ was stated in \cite[Propositions 2.3 and 2.6]{PerVill04} only for full tensor norms. In order to obtain the destruction test for s-tensor norms we need symmetric versions of \cite[Propositions 2.3 and 2.6]{PerVill04} (which are of independent interest).
Floret in \cite{Flo01(extension)} showed that for every s-tensor norm $\beta$ of order $n$ there exist a full tensor norm $\Phi(\beta)$ of order $n$ which is equivalent to $\beta$ when restricted on symmetric tensor products (i.e. there is a constant $d_n$ depending only on $n$ such that $d_n^{-1} \Phi(\beta)|_s  \leq \beta \leq d_n \Phi(\beta)|_s$ in $\otimes^{n,s}E$ for every Banach space $E$). As a consequence a large part of the isomorphic theory of norms on symmetric tensor products can be deduced from the theory of ``full'' tensor norms, which usually is easier to handle.
Using Floret's techniques it is now easy to obtain:

\begin{theorem} \label{Perez García- Villanueva para el caso simétrico} Let $\beta$ be a s-tensor norm of order $n$. If $\big( \widetilde{\otimes}^{n,s} c_0, \beta \big)$ has unconditional basis, then $\beta$ has to be equivalent to $\varepsilon_s$ on $\otimes^{n,s}c_0$. If  $\big( \widetilde{\otimes}^{n,s} \ell_1, \beta \big)$ has unconditional basis, then $\beta$ has to be equivalent to $\pi_s$ on $\otimes^{n,s} \ell_1$.
\end{theorem}

\begin{proof}
First notice that $\big( \widetilde{\otimes}^{n,s} \ell_2^n(c_0), \beta \big) \simeq \big( \widetilde{\otimes}^{n,s} \ell_\infty^n(c_0), \beta \big) \simeq \big( \widetilde{\otimes}^{n,s} c_0, \beta)$ has the Gordon-Lewis property.
Since $\big( \widetilde{\otimes}^n c_0, \Phi(\beta) \big)$ is a complemented subspace of $\big( \widetilde{\otimes}^{n,s} \ell_2^n(c_0), \beta \big)$ we have, by Theorem~\ref{symm-Gl-Incond}, that $\big( \widetilde{\otimes}^n c_0, \Phi(\beta) \big)$ has unconditional basis. Thanks to Pérez-García and Villanueva's result \cite[Proposition 2.3]{PerVill04} we can conclude that $\Phi(\beta) \sim \varepsilon$. Now using the fact that $\Phi(\beta)|_s \sim \beta$ \cite[Theorem 2.3.]{Flo01(extension)} and $\varepsilon|_s \sim  \varepsilon_s$ we get $\beta \sim \varepsilon_s$.

With an analogous proof we obtain that if $\big( \widetilde{\otimes}^{n,s} \ell_1, \beta \big)$ has unconditional basis, $\beta$ must be equal (up to constants) to $\pi_s$.
\end{proof}

%
%
%
%
%
\bigskip
Using Proposition~\ref{Perez García- Villanueva para el caso simétrico}, \cite[Theorem 2.6.]{PerVill05(unconditional)}, Theorem~\ref{symm-Gl-Incond},  and proceeding in a similar way as in Theorem~\ref{Destruction test} we have the following analogous to Theorem~\ref{Destruction test}.

\begin{theorem} \textbf{Destruction Test (symmetric version):} \label{Symmetric Destruction test}
Let $\beta$ be an s-tensor norm of order $n$.
The tensor product $\big(\widetilde{\otimes}^{n,s} E, \beta \big)$ does not have unconditional basis for any Banach space $E$ with unconditional basis if and only if $\beta$ is not equivalent to $\varepsilon_s$, $\sigma_2$ and $\pi_s$ on $\otimes^{n,s} c_0$, $\otimes^{n,s} \ell_2$ and $\otimes^{n,s} \ell_1$ respectively.
\end{theorem}

Let $\beta$ be an s-tensor norm, combining both versions of the destruction test (or using Floret's construction again) we have that $\Phi(\beta)$ destroys unconditionality if and only if $\beta$ destroys unconditionality. Finally, since $/ \Phi(\beta) \s \sim \Phi(/ \beta \s)$ (Floret's construction preserves injective hulls, see \cite{CaGa10}), we can use the previous comment and duality to give the symmetric version of Theorem~\ref{inyectivas y proyectivas destruyen}:
\begin{theorem}
Nontrivial injective and nontrivial projective s-tensor norms destroy unconditionality.
\end{theorem}
\bigskip

We finish this section with the proof of Lemma~\ref{GL-chevet}. We follow
the procedure of \cite[Theorem~3]{DeDiGaMa01}).

\begin{proof}(of Lemma~\ref{GL-chevet})

For the lower estimate, notice first that $\ell_1 \otimes_{\varepsilon} \ell_2 \simeq \ell_1 \otimes_{\alpha} \ell_2$  since $\ell_1$ and $\ell_2$ have cotype 2 \cite[Exersice 31.2]{DefFlo93}. Then, $\ell_1^m \otimes_{/ \pi \s} \ell_2^m$ is isomorphic to $\ell_1^m \otimes_{\varepsilon} \ell_2^m$ with constants independent of $m$.
So we have to estimate $gl \big( \ell_1^n \otimes_{\varepsilon} \ell_2^n \big)$.

Observe that
\begin{equation*}
m^{3/2} \leq  \| \sum_{i,j}^m  e_i \otimes e_j \|_{\ell_1^m \otimes_{\varepsilon} \ell_2^m}.
\end{equation*}
As a matter of fact,
$$\| \sum_{i,j}^m  e_i \otimes e_j \|_{\ell_1^m \otimes_{\varepsilon} \ell_2^m} = \sup_{a \in B_{\ell_\infty^m}, b \in B_{\ell_2^m} } \big| \sum_{i,j}^m a_i b_j \big| \geq m \ \sup_{ b \in B_{\ell_2^m}} \big| \sum_{j}^m  b_j \big| = m m^{1/2}= m^{3/2}.$$

We now consider the aleatory matrices
\begin{align*}
R : \Omega \rightarrow \ell_1^m \otimes_{\varepsilon} \ell_2^m \; \; & R(\omega) := \sum_{i,j}^m r_{i,j}(\omega) e_i \otimes e_j \\
G: \Omega \rightarrow \ell_1^m \otimes_{\varepsilon} \ell_2^m \; \; & G(\omega) := \sum_{i,j}^m g_{i,j}(\omega) e_i \otimes e_j, \\
\end{align*}
where $(\Omega,\mu)$ is a probability space and $r_{i,j}$'s and $g_{i,j}$'s forms a family of $m^2$ Bernoulli and Gaussian variables on $\Omega$, respectively.

Then, for all $\omega \in \Omega$
\begin{align*}
m^{3/2} \leq  \| \sum_{i,j}^m  e_i \otimes e_j \|_{\ell_1^m \otimes_{\varepsilon} \ell_2^m} & = \| \sum_{i,j}^m  r_{i,j}(\omega) r_{i,j}(\omega) e_i \otimes e_j \|_{\ell_1^m \otimes_{\varepsilon} \ell_2^m} \\
& \leq \chi((e_i \otimes e_j)_{i,j}) \|R(\omega)\|_{\ell_1^m \otimes_{\varepsilon} \ell_2^m}.
\end{align*}

On the other hand, we know that  $\chi((e_i \otimes e_j)_{i,j}) \leq 2^{3} gl( \ell_1^n \otimes_{\varepsilon} \ell_2^n)$ by Theorem~\ref{full-Gl_incond}.
Therefore, for every $\omega \in \Omega$ we have:
$$m^{3/2} \leq  \| \sum_{i,j}^m  e_i \otimes e_j \|_{\ell_1^m \otimes_{\varepsilon} \ell_2^m} \leq 2^{3} gl( \ell_1^m \otimes_{\varepsilon} \ell_2^m) \|R(\omega)\|_{\ell_1^m \otimes_{\varepsilon} \ell_2^m}.
$$
Integrating the last expression,
$$m^{3/2} \leq  \| \sum_{i,j}^m  e_i \otimes e_j \|_{\ell_1^m \otimes_{\varepsilon} \ell_2^m} \leq 2^{3} gl( \ell_1^m \otimes_{\varepsilon} \ell_2^m) \int_{\Omega}  \|R(\omega)\|_{\ell_1^m \otimes_{\varepsilon} \ell_2^m} \; d\mu.
$$
Now, since Gaussian averages $L=1/{M}_1$ dominate, up to a uniform constant, Bernoulli averages \cite[Page 15.]{Tomczak-Jae}, \cite[Proposition 12.11.]{DieJarTon95} we get:
$$ \int_{\Omega} \|R(\omega)\|_{\ell_1^m \otimes_{\varepsilon} \ell_2^m} \; d\mu \leq L \int_{\Omega} \|G(\omega)\|_{\ell_1^m \otimes_{\varepsilon} \ell_2^m} \; d\mu.$$

It is time to use Chevet inequality:

\noindent\textbf{Chevet Inequality} \cite[(43.2)]{Tomczak-Jae}:
Let $E$ and $F$ be Banach spaces. Fix $x_1', \dots, x_m' \in E'$ and $y_1, \dots, y_m \in F$.
If $\{g_{ij}\},\{g_{i}\},\{g_{j}\}$ are independent Gaussian random variables in some probability space $(\Omega,\mu)$.
Then,
\begin{align*}
\int_{\Omega} \|\sum_{i,j=1}^m g_{ij} x_i' \otimes y_j\| \;d\mu & \leq  b \; \mbox{sup}_{{ \|x\|}\leq  1} \left( \sum_{i=1}^m |x_i'(x)|^2\right)^{1/2} \int_{\Omega} \| \sum_{j =1}^m g_{j} y_j \| \;d\mu \; \\
& + b \; \mbox{sup}_{{ \|y'\|}\leq  1} \left( \sum_{j=1}^m |y'(y_j)|^2\right)^{1/2} \int_{\Omega} \|\sum_{i =1}^m g_{i} x_i' \| \;d\mu.
\end{align*}
To conclude with our estimations it remains to observe that $$\|\sum_{i,j}^m g_{i,j}(\omega) e_i \otimes e_j\|_{\ell_1^m \otimes_{\varepsilon} \ell_2^m} = \| \sum_{i,j}^m g_{i,j}(\omega) e_i \otimes e_j \|_{ L( \ell_{\infty}^m,\ell_2^m) }.$$
Then,
\begin{eqnarray*}\int_{\Omega} \|G(\omega)\|_{\ell_1^m \otimes_{\varepsilon} \ell_2^m} \; d\mu & \leq b \sup_{x \in B_{\ell_\infty^m}} \big( \sum_{i=1}^m |x_i|^2 \big)^{1/2} \int_{\Omega} \| \sum_{j=1}^m g_j(\omega) e_j \|_{\ell_2^m} \; d\mu(\omega)\\ & + b \sup_{y \in B_{\ell_2^m}} \big( \sum_{j=1}^m |y_j|^2 \big)^{1/2} \int_{\Omega} \| \sum_{i=1}^m g_i(\omega) e_i \|_{\ell_1^m}.\end{eqnarray*}
Using \cite[Proposition 45.1]{Tomczak-Jae} we have that the last member is less or equal to $Cn$, where $C$ is a constant.
We have shown that $m^{1/2} \leq C gl \big( \ell_1^m \otimes_{\varepsilon} \ell_2^m \big)$.

The upper estimate follows from the fact that $d(\ell_1^m \otimes_{\varepsilon} \ell_2^m,\ell_1^m \otimes_{\varepsilon} \ell_{\infty}^m)  \leq d(\ell_2^m,\ell_{\infty}^m)= m^{1/2}$ together with $\chi(\ell_1^m \otimes_{\varepsilon} \ell_{\infty}^m)=1$ (since $\chi(\ell_1^m)=1$).

\end{proof}

\section{Unconditionality in ideals of polynomials and multilinear forms}\label{section-ideals}

Let $\mathcal{P}^n(E)$ denote the space of all continuous $n$-homogeneous scalar-valued polynomials on $E$, endowed with ``supremum on the ball'' norm.

Let us recall the definition of polynomial ideals
\cite{Flo02(On-ideals)}. A \emph{Banach ideal of continuous scalar valued
$n$-homogeneous polynomials} is a pair
$(\mathcal{Q}^n,\|\cdot\|_{\mathcal Q^n})$ such that:
\begin{enumerate}
\item[(i)] $\mathcal{Q}^n(E)=\mathcal Q^n \cap \mathcal{
P}^n(E)$ is a linear subspace of $\mathcal{P}^n(E)$ and $\|\cdot\|_{\mathcal Q^n}$ is a norm which makes the pair
$(\mathcal{Q}^n,\|\cdot\|_{\mathcal Q^n})$ a Banach space.

\item[(ii)] If $T\in \mathcal{L} (E_1,E)$, $P \in \mathcal{Q}^n(E)$ then $P\circ T\in \mathcal{Q}^n(E_1)$ and $$ \|
P\circ T\|_{\mathcal{Q}^n(E_1)}\le  \|P\|_{\mathcal{Q}^n(E)} \| T\|^n.$$

\item[(iii)] $z\mapsto z^n$ belongs to $\mathcal{Q}^n(\mathbb K)$
and has norm 1.
\end{enumerate}

\bigskip

The classes $\mathcal{P}_{N}^n$ and $\mathcal{P}_{I}^n$  of nuclear and integral polynomials are examples of Banach ideals of polynomials (see \cite{Flo97} for definitions).

We say that an $n$-homogeneous polynomial $P : \ell_2 \to \mathbb{K}$ is Hilbert Schmidt if
$$ \big( \sum_{k_1, \dots, k_n =1}^\infty |\check{P}(e_{k_1},\dots,e_{k_n})|^2 \big)^{1/2} < \infty,$$
where $\check{P}$ is the symmetric $n$-linear form associated to $P$.
The space of all such polynomial will be denoted by $\mathcal{P}_{HS}^n(\ell_2)$ with the norm $\|P\|_{\mathcal{P}_{HS}^n(\ell_2)} =  \big( \sum_{k_1, \dots, k_n =1}^\infty |\check{P}(e_{k_1},\dots,e_{k_n})|^2 \big)^{1/2}$.

We begin with a reformulation of the Destruction test in terms of ideals of polynomials:

\begin{proposition} \label{ideal destruction test} If $\mathcal Q^n$ is a Banach ideal of $n$-homogeneous polynomials, the following are equivalent:

(i) For any Banach space $E$ with unconditional basis,  $\mathcal{Q}^n(E)$ fails to have the Gordon-Lewis property.

(ii) $\|Id:\mathcal{Q}^n(\ell_\infty^m)\to \mathcal{P}_{I}^n(\ell_\infty^m)\|\to\infty$,  $\|Id:  \mathcal{P}^n(\ell_1^m)    \to \mathcal{Q}^n(\ell_1^m) \|\to\infty$ and \\ $\max \big(\|Id:\mathcal{Q}^n(\ell^m_2)\to \mathcal{P}_{HS}^n(\ell^m_2)\|,
\|Id:\mathcal{P}_{HS}^n(\ell^m_2) \to  \mathcal{Q}^n(\ell^m_2)\| \big) \to\infty$ as $m\to\infty$.

If $\mathcal{Q}$ is maximal, this is also equivalent to

(iii)  $\mathcal{Q}^n(\ell_1) \neq \mathcal{P}^n(\ell_1)$, $\mathcal{Q}^n(\ell_2) \neq \mathcal{P}_{HS}^n(\ell_2)$ and $\mathcal{Q}^n(c_0) \neq \mathcal{P}_{I}^n(c_0)$.
\end{proposition}

\begin{proof} It is clear that (i) implies any of the other statements.
To see that (ii) implies (i), by Tzafriri's result (Theorem~\ref{tzafriri}) it is enough to see that $gl\big( \mathcal{Q}^n(\ell_p^m) \big) \to \infty$ as $m \to \infty$ for all these $p=1,2,\infty$.
We can suppose $p=1$, the other cases being completely analogous. Let $\beta$ be the s-tensor norm associated to $\mathcal{Q}^n$. Since $\mathcal{Q}^n(\ell_\infty^m) = \big( \otimes^{n,s} \ell_1^m, \beta \big)$, if $gl \Big(\big( \otimes^{n,s} \ell_1^m, \beta \big) \Big)$ were uniformly bounded we would have that $\beta \sim \pi_s$ on $\otimes^{n,s} \ell_1$ by Theorem~\ref{Perez García- Villanueva para el caso simétrico}. Therefore, the norms of $\mathcal{Q}^n(\ell_\infty^m)$ and $\mathcal{P}_{I}^n(\ell_{\infty}^m)$ would be equivalent (with constants independent of $m$), a contradiction.

The maximality of $\mathcal Q^n$ assures that, if $\|Id:\mathcal{P}^n(\ell_1^m)\to \mathcal{Q}^n(\ell_1^m)\|$ is uniformly bounded on $m$, then $\mathcal{P}^n(\ell_1) = \mathcal{Q}^n(\ell_1)$ and, of course, the converse is also true. The same holds for the other two conditions in (ii) and (iii).
\end{proof}


For a Banach space $E$ with unconditional basis $(e_j)_{j=1}^\infty$, the authors of \cite{DefKal05} studied when $\mathcal{P}^n(E)$ was isomorphic to a Banach lattice. It turned out that this happens precisely when the monomials $e_{\alpha}'$ form an unconditional basic sequence.

The same holds for a maximal Banach polynomial ideals:

\begin{proposition}\label{monomios-incond}
Let $\mathcal{Q}^n$ be a maximal ideal of $n$-homogeneous polynomials and $E$ be a Banach space with unconditional basis $(e_j)_{j=1}^\infty$. The following are equivalent:
\begin{enumerate}
\item The monomials $e_{\alpha}'$ form an unconditional basic sequence in $\mathcal{Q}^n(E)$.
\item $\mathcal{Q}^n(E)$ is isomorphic to a Banach lattice.
\item $\mathcal{Q}^n(E)$ has the Gordon-Lewis property.
\end{enumerate}
\end{proposition}

The proposition can be proved similarly to  \cite[Proposition 4.1]{DefKal05} and with the help of the following lemma, which is a consequence of the polynomial version of the Density Lemma \cite{DefFlo93}:

\begin{lemma}\label{estimacion de la norma con truncados}
Let $\mathcal{Q}^n$ be a maximal ideal of $n$-homogeneous polynomials and $E$ a Banach space with monotone basis $(e_j)_{j=1}^\infty$.
A continuous $n$-homogeneous polynomial $Q$ belongs to $\mathcal{Q}^n(E)$ if and only if  $\sup_{k \in \mathbb{N}} \|Q \circ P_k \|_{Q^n(E)} < \infty$, where $P_k$ stands for the canonical projection to $E_k:= [e_j : 1 \leq j \leq k]$.

In that case,
$$\|Q\|_{\mathcal{Q}^n(E)} = \sup_{k \in \mathbb{N}} \|Q \circ P_k \|_{\mathcal{Q}^n(E)} = \sup_{k \in \mathbb{N}} \|Q{|_{E_k}} \|_{\mathcal{Q}^n(E_k)}.$$
\end{lemma}

\bigskip
Now we present some examples of Banach polynomial ideals that destroy the Gordon-Lewis property (in the sense of the Proposition~\ref{ideal destruction test}). An immediate consequence of Theorem~\ref{s-tensor norms naturales destruyen} is the following:

\begin{proposition}If $\mathcal{Q}^n$ is a Banach ideal of $n$-homogeneous polynomials associated to a nontrivial injective norm or a nontrivial projective norm, then $\mathcal{Q}^n(E)$ does not have the Gordon-Lewis property for any Banach space with unconditional basis (equivalently the monomial basic sequence is not unconditional). This holds, in particular, for any ideal with a nontrivial natural associated s-tensor norm.
\end{proposition}

As an example of the latter, we have:

\textbf{The ideal of extendible polynomials:} A polynomial $P \in \mathcal{P}^n$ is \textbf{extendible} \cite{KirRy98}
if for any Banach space $G$ containing $E$ there exists
$\widetilde P\in {\mathcal P}^n(G)$ an extension of $P$. We will
denote the space of all such polynomials by ${\mathcal
P}_e^n(E)$. For $P\in {\mathcal P}_e^n(E)$, its extendible
norm is given by
$$\begin{array}{rl}
\Vert P\Vert _{{\mathcal P}_e^n(E)}=\inf \{C>0:  & \mbox{for
all }G\supset E
\mbox{ there is an extension of }P\mbox{ to }G  \\
& \mbox{ with norm}\leq C\}.
\end{array}$$

It is known that $\mathcal{P}^n_e(E) = \big( \otimes^{n,s} E, {/ \pi_s \s} \big)'$ (see \cite[Section 20.7]{DefFlo93} and \cite{Carando99}), so we have:

\begin{example}\label{ejemplo-extendibles}
If $E$ is a Banach space with unconditional basis then $\mathcal{P}_e^n(E)$ does not have the Gordon-Lewis property. Therefore, by Proposition~\ref{monomios-incond}, the monomial basic sequence is not unconditional in $\mathcal{P}_e^n(E)$.
\end{example}

\textbf{The ideal of $r$-dominated polynomials:}
For $x_1,\dots,x_m\in E$, we define
$$
w_{r} \big((x_{i})_{i=1}^m\big) = \sup_{x' \in B_{E'}} \left(
\sum_{i} |<x', x_i>|^{r} \right)^{1/r}.
$$
A polynomial $P \in P^n(E)$ is
\textbf{$r$-dominated} (for $r \geq n$) if there exists $C>0$ such that for every
finite sequence $(x^{i})_{i=1}^{m} \subset E$ the following holds
\[
\left( \sum_{i=1}^{m} |P (x_i)|^{\frac{r}{n}}
\right)^{\frac{n}{r}} \leq C w_{r}((x_i)_{i=1}^{m})^n.
\]
We will denote the space of all such polynomials by  $\mathcal{D}_{r}^n(E)$.
The least of such constants $C$ is called the
\emph{$r$-dominated norm} and denoted $\|P\|_{\mathcal{D}_r ^n(E)}$.

\begin{example} \label{polinomios r dominados}
If $E$ is a Banach space with unconditional basis and $r \geq n$, then $\mathcal{D}_r^n(E)$ does not have the Gordon-Lewis property and the monomial basic sequence is never unconditional.
\end{example}

\begin{proof}
By Proposition~\ref{ideal destruction test} we must show that $\mathcal{D}_r^n(\ell_1) \neq \mathcal{P}^n( \ell_1)$, $\mathcal{D}_r^n(\ell_2) \neq \mathcal{P}_{HS}^n(\ell_2)$, $\mathcal{D}_r^n(c_0) \neq \mathcal{P}_{I}^n( c_0)$.

If $\mathcal{D}_r^n(\ell_1) = \mathcal{P}^n( \ell_1)$, using \cite[Lemma 1.5]{CarDimMurN} we would have that $\mathcal{D}_r^2(\ell_1) = \mathcal{P}^2( \ell_1)$ (since $\mathcal{D}_r$ and $\mathcal{P}$ are coherent sequences of polynomial ideals \cite[Examples 1.9, 1.13]{CarDimMurN}). By \cite[Proposition 12.8]{DefFlo93} we would have: $\mathcal{D}_r^2(\ell_1)=\mathcal{D}_2^2(\ell_1)=\mathcal{P}_e^2( \ell_1)$, but we already know that $\mathcal{P}_e^2( \ell_1)$ cannot be equal to $\mathcal{P}^2(\ell_1)$.

Using coherence again, it is easy to show that $\mathcal{P}_{HS}^n(\ell_2) \not\subset \mathcal{D}_r(^n\ell_2)$ (recall that Hilbert Schmidt polynomials coincide with multiple 1-summing polynomials, which form a coherent sequence of ideals \cite[Example 1.14]{CarDimMurN}): if $ \mathcal{P}_{HS}^n(\ell_2) \subset \mathcal{D}_r^n(\ell_2)$, we would have $\mathcal{P}_{HS}(^2 \ell_2)  \subset   \mathcal{D}_r(^2\ell_1)  =\mathcal{D}_2^2(\ell_2)=    \mathcal{P}_e^2( \ell_2)$, which is not true, for example, by (\ref{eq eta sigma 2}) and duality.

Similarly, $\mathcal{D}_r^n( c_0) \neq  \mathcal{P}_{I}^n( c_0)$ ($\mathcal{P}_{I}$ is also a coherent sequence \cite[Example 1.11]{CarDimMurN}).
\end{proof}

\textbf{The ideal of $r$-integral polynomials:}
If $\mu$ is a finite, positive measure on $\Omega$ and $n \leq r \leq \infty$, the $n$-th integrating polynomial $q_{\mu,r}^n (f) : = \int f^n d\mu$. It is straightforward to see that $\|q_{\mu,r}^n\| = \mu(\Omega)^{1/s}$ where $s =(\frac{r}{n})'$. A polynomial $P \in \mathcal{P}^n(E)$ is \textbf{$r$-integral} \cite{Flo02(On-ideals)} if it admits a factorization
$$ P : E \overset{T} \longrightarrow L_r ( \Omega) \overset{q_{\mu,r}^n} \longrightarrow \mathbb{K}$$
with a finite, positive measure $\mu$ and $T \in \mathcal{L}(E, L_r ( \Omega))$. We will
denote the space of all such polynomials by $\mathcal{I}_r^n(E)$.
And
$$\|P\|_{\mathcal{I}_r^n(E)}= \inf \{ \|T\|^n \|q_{\mu,r}^n\| \; : \; P = q_{\mu,r}^n \circ T \; \mbox{as before}\}.$$

\begin{example}
If $E$ be a Banach space with unconditional basis and $r \geq n$, then  $\mathcal{I}_r^n(E)$ does not have the Gordon-Lewis property.
\end{example}

\begin{proof}
As in the proof of \cite[Theorem 3.5]{CarDimSev07} we can see that, if $M$ is a finite dimensional space, then $(\mathcal{D}_r^n)^{*}(M)= \mathcal{I}_r^n(M)$. Let $p \in \{ 1,2, \infty \}$, we have that $gl\big(\mathcal{I}_r^n(\ell_p^m)\big)= gl\big(\mathcal{D}_r^n( \ell_p^m)\big)$, which we already know by the previous example that this goes to $\infty$ with $m$.
\end{proof}

Note that in the proofs of the previous examples  we have actually shown the following limits, which we will use below: \begin{eqnarray}\label{id-HS-Dr} \|Id:\mathcal{P}_{HS}^n(\ell_2^m)\to \mathcal{D}_r^n(\ell_2^m)\|\to\infty \\ \label{id-Ir-HS} \|Id: \mathcal{I}_r^n(\ell_2^m) \to \mathcal{P}_{HS}^n(\ell_2^m)\|\to\infty\end{eqnarray} as $m$ goes to infinity.

\bigskip

\textbf{The ideal of  $r$-factorable polynomials:}
For $n \leq r \leq \infty$, a polynomial $P \in \mathcal{P}^n(E)$ is called \textbf{$r$-factorable} \cite{Flo02(On-ideals)} if there is a positive measure space $(\Omega, \mu)$, an operator $T \in \mathcal{L}\big(E, L_r(\mu)\big)$ and $Q \in \mathcal{P}^n \big( L_r(\mu) \big)$ with $P = Q \circ T$. The space of all such polynomials will be denoted by $\mathcal{L}_r^n(E)$.
With
$$ \|P\|_{\mathcal{L}_r^n(E)} = \inf \{ \|Q\| \|T\|^n \; : \; P: E \overset T \longrightarrow L_r(\mu) \overset Q \longrightarrow \mathbb{K}\}.$$

\begin{example}
Let $E$ be a Banach space with unconditional basis and $r \geq n$, then $\mathcal{L}_r^n(E)$ does not have the Gordon-Lewis property and the monomial basic sequence is not unconditional in $\mathcal{L}_r^n(E)$.
\end{example}

\begin{proof}
By \cite[Theorem 3.5]{CarDimSev07} and then \cite[Proposition 4.3.]{Flo02(On-ideals)}, we have $\mathcal{D}_r^{*}=\mathcal{I}_r^{max} \subset \mathcal{L}_r$ ($\mathcal{L}_r$ is maximal \cite[Proposition 3.1]{Flo02(On-ideals)}). Therefore, using Proposition~\ref{ideal destruction test} and Equation~(\ref{id-Ir-HS}), we have $\|Id:\mathcal{L}_r^n(\ell_\infty^m)\to \mathcal{P}_{I}^n(\ell_\infty^m)\|\to\infty$ and  $\|Id: \mathcal{L}_r^n(\ell_2^m) \to \mathcal{P}_{HS}^n(\ell_2^m)\|\to\infty$.
It remains to show that $\mathcal{L}_r^n( \ell_1) \neq \mathcal{P}^n( \ell_1)$.
We will show this first for $n=2$. Suppose this happens, then every symmetric operator $T: \ell_1 \to \ell_\infty$ would factorize by a reflexive Banach Space, then must be weakly compact, a contradiction to the fact that $\ell_1$ is not symmetrically Arens regular \cite[Section 8]{AroColGam91}.
For $n \geq 3$ we use coherence for composition ideals \cite[Proposition 3.3]{CarDimMurN} since $\mathcal{L}_r  = \mathcal{P} \circ \Gamma_r$ \cite[3.5.]{Flo97}.
\end{proof}

In \cite{DefKal05}, Defant and Kalton showed that the space $\mathcal P^n(E)$ of all $n$-homogeneous polynomials cannot have unconditional basis whenever $E$ is a Banach space with unconditional basis. However, $\mathcal P^n(E)$ can have the Gordon-Lewis property (for example, when $E=\ell_1$). When this happens, $\mathcal P^n(E)$ is not separable and therefore it has no basis at all. One may wonder if there are other ideals with that property: that never have unconditional basis but sometimes enjoy the Gordon-Lewis property. We will present such an example but first we extend the range of ideals for which  \cite[Proposition 3.2.]{DefKal05} apply. For each $m$, we define $P_m\in\mathcal P(^n\ell_2)$ by $$P_m(x)=\sum_{j=1}^m x_j^n.$$

\begin{proposition} \label{ideal separable implica contencion}
Let $\mathcal{Q}^n$ a polynomial ideal such that $(\|P_m\|_{\mathcal{Q}^n(\ell_2)})_m$ is uniformly bounded.
If $(\mathcal{Q}^n)^{max}(E)$ is separable, then $E$ must contain  $(\ell_2^m)_{m=1}^\infty$ or $(\ell_{\infty}^m)_{m=1}^\infty$ uniformly complemented.
\end{proposition}

\begin{proof}
Let $(e_k)_{k=1}^\infty$ be an unconditional basis of $E$. By the proof of \cite[Proposition 3.2.]{DefKal05} we know that if $E$ does not contain any of the sequences $(\ell_2^m)_{m=1}^\infty$, $(\ell_{\infty}^m)_{m=1}^\infty$ uniformly complemented, we may extract a subsequence $(f_{j})_{j=1}^\infty$ of $(e_k)_{k=1}^\infty$ such that for any $x \in F:=\overline{[(f_j)]}$ we have $\sum_{j=1} |f_j'(x)|^2 < \infty $, where $(f_j')_j$ is the corresponding subsequence of the dual basic sequence. This means that, as sequence spaces, we have a continuous  inclusion  $i: F  \hookrightarrow \ell_2$. For $x \in F$, we define $Q_m(x)=\sum_j^m f_j'(x)^n$. We have
 $$\|Q_m \|_{\mathcal{Q}^n(F)} = \|P_m\circ i \|_{\mathcal{Q}^n(F)} \le \|P_m \|_{\mathcal{P}^n(\ell_2)}\|i\|^n,$$which is bounded uniformly on $m$.
It follows from \cite[Lemma 5.4]{CaDimSe08} that $(\mathcal{Q}^n)^{max}(F)$ cannot be separable. Therefore, $(\mathcal{Q}^n)^{max}(E)$ cannot be separable either since $F$ is a complemented subspace of $E$.
\end{proof}

The uniform bound for $(\|P_m\|_{\mathcal{Q}^n(\ell_2)})_m$ is necessary for the result to be true, as the following example shows.

\begin{example}
Let $E$ be a reflexive Banach space with unconditional basis. Since $ \s \varepsilon /$ is a projective 2-fold tensor norm with the Radon-Nikod\'ym property \cite[Theorem 33.5]{DefFlo93}, we have $$E' \otimes_{\s \varepsilon /} E' = (E \otimes_{/ \pi \s} E)' = \mathcal L_e^2(E).$$ Therefore, $\mathcal L_e^2(E)$ and $\mathcal P_e^2(E)$ have Schauder bases \cite{GreRy05} (which we already know are not unconditional) and are consequently separable. If we take $E$ to be the dual of the original Tsirelson's space, $E$ does not contain either $\ell_2^m$ nor $\ell_{\infty}^m$ uniformly complemented \cite[Pages 33 and 66]{Casazza-Shura(Tsirelson)}.
\end{example}

\begin{corollary} \label{Q max no tiene base incondicional}
Let $\mathcal{Q}^n$ be a maximal Banach ideal of $n$-homogeneous polynomials such that $(\|P_m\|_{\mathcal{Q}^n(\ell_2)})_m$ is uniformly bounded. Suppose also that not ever polynomial in $\mathcal Q^n(c_0)$ is integral. If $E$ or its dual has unconditional basis, then $\mathcal{Q}^n(E)$ does not have unconditional basis.
\end{corollary}

\begin{proof} Suppose first that $E$ has unconditional basis. If $\mathcal{Q}^n(E)$ is separable, by Proposition~\ref{ideal separable implica contencion} $E$ must contain either $(\ell_{\infty}^m)_{m=1}^\infty$ or $(\ell_2^m)_{m=1}^\infty$ uniformly complemented.
If $E$ contains the sequence $(\ell_{\infty}^m)_{m=1}^\infty$ uniformly complemented, since not every polynomial on $c_0$ is integral, we have by the proof of Proposition~\ref{ideal destruction test} $gl(\mathcal Q^n(\ell_\infty^m))\to\infty$ as $m\to \infty$, so $\mathcal{Q}^n(E)$ cannot have the Gordon-Lewis property. If $E$ contains $(\ell_2^m)_{m=1}^\infty$ uniformly, since $(\|P_m\|_{\mathcal{Q}^n(\ell_2)})_m$ is uniformly bounded and $(\|P_m\|_{\mathcal{P}_{HS}^n(\ell_2)})_m=\sqrt m$, we can conclude that $gl(\mathcal Q^n(\ell_2^m))\to\infty$ as $m\to\infty$. Therefore, if $E$ is reflexive, $\mathcal{Q}^n(E)$ either fails the Gordon-Lewis property or is non-separable. In any case, it has no unconditional basis.

If $E'$ has unconditional basis and is reflexive, then $E$ also has unconditional basis and we are in the previous case. If $E'$ is not reflexive and has unconditional basis, then $E'$ contains complemented copies of $c_0$ or $\ell_1$. If it contains $c_0$, it also contains $\ell_\infty$, so $E'$ is not separable, and neither is $\mathcal{Q}^n(E)$. If $E'$ contains $\ell_1$ and we denote by $\beta$ the s-tensor norm associated to $\mathcal Q$, we have that $\mathcal Q(E)$ contains the spaces $\otimes_{\beta}^{n,s}\ell^m_1$ which are uniformly isomorphic to $\mathcal Q^n(\ell_\infty^m)$. As in the reflexive case, the Gordon-Lewis constant of $\mathcal Q^n(\ell_\infty^m)$ goes to infinity, so $\mathcal Q^n(E)$ does not have the Gordon-Lewis property.
\end{proof}

As a consequence of the previous corollary, we have that $\mathcal P^n(E)$ cannot have an unconditional basis for any Banach space $E$ that has (or its dual has) unconditional basis. Since there are Banach spaces without unconditional basis whose duals have one (see for example the remark after \cite[1.c.12.]{Lind-Tza(Classical-VOL_I)}), this somehow extends their result and answers Dineen's question as it was originally posed. However, it should be stressed that our arguments strongly rely on Defant and Kalton's work.

Another consequence is the following: suppose that $E'$ has a Schauder basis $(e_j')_{j=1}^\infty$ and $\mathcal Q^n$ is as in the previous corollary. Then, the monomials associated to $(e_j')_{j=1}^\infty$ cannot be an unconditional basis of $\mathcal Q^n(E)$. Indeed, should the monomials be an unconditional sequence, then $(e_j')_{j=1}^\infty$ would be also unconditional, so we can apply Corollary~\ref{Q max no tiene base incondicional}.

\medskip

Now we present another example of a Banach ideal of polynomials which can have the Gordon-Lewis property but that never has unconditional basis, just as $\mathcal P^n$:

Consider $\mathcal{Q}^n$ the ideal given by $\mathcal{Q}^n = D_n^n \circ \Gamma_{\infty}^{-1}$ (here we follow the notation of \cite[25.6]{DefFlo93} for quotient ideals). More precisely, a polynomial $P$ belongs to $\mathcal{Q}^n(E)$ if there exists a constant $C>0$ such that for every $\infty$-factorable operator $T: F \to E$  with $\gamma_{\infty}(T) \leq 1$, the composition $P \circ T$ is $n$-dominated and $\| P \circ T \|_{D_n^n} \leq C$.
We define $$\|P\|_{\mathcal{Q}^n} := \sup \{ \| P \circ T \|_{D_n^n} : \gamma_\infty(T) \leq 1 \},$$ where $D_n^n$ is the ideal of $n$-dominated polynomials.

It is not hard to see that $\mathcal{Q}^n$ is in fact a Banach ideal of of $n$-homogeneous polynomials. Also, we have that $\mathcal{Q}^n(\ell_1) = \mathcal{P}^n(\ell_1)$. Indeed,
take $P \in \mathcal{P}^n(\ell_1)$ and $T \in \Gamma_{\infty}(F, \ell_1)$ with unit norm and let us find a constant $C$ such that $\| P \circ T \|_{D_n^n} \leq C$.
If $S : F \to L_{\infty}(\mu)$ and $R : L_{\infty}(\mu) \to \ell_1$ are operators which satisfy $\|S\| \|R\| \leq 2$ and $T = S \circ R$,
then $P \circ T = P \circ R \circ S $.
By Grothendieck theorem, $R$ is $n$-summing and
$\pi_n(R) \leq K_G \|R\|.$ Since $D_n^n$ is the composition ideal $\mathcal{P}^n \circ \Pi_n$ \cite{Sch91} we have that $\|P \circ R \|_{D_n^n} \leq  K_G^n \|P\| \|R\|^n$. Therefore
 $\|P \circ T\|_{D_n^n} \leq K_G^n \|P\| \|R\|^n \|S\|^n \leq (2 K_G)^n \|P\|$
and we are done.

Using a similar argument it can be shown that $\mathcal{Q}^n(\ell_2)= \mathcal{P}^n( \ell_2)$, so the sequence $(\|P_m\|_{\mathcal{Q}^n(\ell_2)})_m$ is uniformly bounded.
We also have that $\mathcal{Q}^n(c_0) = D_n^n ( c_0)\not\supset\mathcal{P}_{I}^n(c_0)$.

Finally, we see that $\mathcal{Q}^n$ is maximal.
Take $P \in (\mathcal{Q}^n)^{max}(E)$ and let us show that $P \in \mathcal{Q}^n(E)$, that is, $\|P \circ T \|_{D_n^n} \leq C$ for every $T \in \Gamma_{\infty}(F, E)$ with $\gamma_\infty(T) \leq 1$.
Since $\mathcal D_n^n$ is a maximal ideal, it is sufficient to prove that
$\|P \circ T|_M \|_{D_n^n} \leq C$ for every $M \in FIN(F)$ and $T$ as before.
But, $P \circ T|_M = P|_{Im (T|_M)} T|_M $ and since $P  \in (\mathcal{Q}^n)^{max}(E)$ we have
$\| P|_N \|_{\mathcal{Q}^n} \leq K$ for every $N \in FIN(E)$.
This means that $\sup_{ \gamma_{\infty}(T) \leq 1} \| P|_N \circ T \|_{D_n^n} \leq K$ and we are done.

Thus, Corollary~\ref{Q max no tiene base incondicional} says that $\mathcal{Q}^n(E)$ has not unconditional basis if $E$ or its dual has unconditional basis. On the other hand, $\mathcal{Q}^n(\ell_1) = \mathcal{P}^n(\ell_1)$ has the Gordon-Lewis property.

\bigskip

We have presented examples of several polynomial ideals that lack the Gordon-Lewis property for any Banach space with unconditional basis. It is easy to obtain the same conclusions for ideals of multilinear forms on a single space. For example, Theorem~\ref{inyectivas y proyectivas destruyen} gives:

\begin{proposition}
Let $\mathcal{U}^n$ be a Banach ideal of $n$-linear forms associated to a nontrivial injective or projective tensor norm. If $E$ has unconditional basis, then $\mathcal{U}^n(E)$ does not have the Gordon-Lewis property.
\end{proposition}

From the previous result and Proposition~\ref{Para mas de 3 arruina} we have:

\begin{example}\label{multi_extend_diferentes}
(i) The space $\mathcal{L}^n_e(E)$ do not have the Gordon-Lewis property for any Banach space $E$ and $n \geq 2$.

(ii) If $E_1, \dots, E_n$ are Banach spaces with unconditional basis ($n \geq 3$) then $\mathcal{L}_e(E_1, \dots, E_n)$ do not have the Gordon-Lewis property.
\end{example}

On the other hand, the comments after Proposition~\ref{Para mas de 3 arruina} shows that we cannot expect (ii) to hold for $n=2$. Moreover, the space $\mathcal{L}_e(c_0,\ell_2)$ not only enjoys the Gordon-Lewis property, in fact it has unconditional basis:  since $\s \varepsilon /$ has the Radon-Nikod\'ym property, $$\mathcal{L}_e(c_0,\ell_2)=(c_0\widetilde{\otimes}_{/\pi\s}\ell_2)'=\ell_1\widetilde{\otimes}_{\s \varepsilon /}\ell_2,$$ and therefore has a monomial basis. Since we have shown that $c_0\otimes_{/\pi\s}\ell_2$ has the Gordon-Lewis property, the monomial basis of $\ell_1 \widetilde{\otimes}_{\s \varepsilon /}\ell_2$ must be unconditional.

\bigskip

An example that does not follow from the injective/projective result is the ideal of \textbf{$r$-dominated multilinear forms:}

Let $r \geq n$, an $n$-linear form $T : E_1 \times \dots \times E_n \to \mathbb{K}$ is \emph{$r$-dominated} if there is a constant $C \geq 0$ such that, however we choose finitely many vector $(x_i^j)_{i=1}^m \in E_j$, we have
$$ \big( \sum_{i=1}^m |T(x_i^1, \dots, x_i^n)|^{r/n} \big)^{n/r} \leq C w_r\big( (x_i^1)_{i=1}^m \big) \dots w_r \big((x_i^n)_{i=1}^m \big).$$
The space of all such $T$ will be denoted $\mathfrak{D}_r^n (E_1, \dots, E_n)$ with the norm $\delta^n_r(T) = \min C$.

Since the ideal of $r$-dominated polynomials $\mathcal{D}_r^n(E)$ is isomorphic to a complemented subspace of $\mathfrak{D}_r^n(E)$, from the polynomial result (Example~\ref{polinomios r dominados}) we obtain:

\begin{example}
Let $E$ be a Banach space with unconditional basis. Then, $\mathfrak{D}_r^n(E):=\mathfrak{D}_r^n(E, \dots, E)$ do not have the Gordon-Lewis property.
\end{example}

Let us mention that for different spaces, we can obtain that dominated multilinear forms behaves exactly as extendible ones in Example~\ref{multi_extend_diferentes}. The case $n=2$ follows from the coincidence between dominated and extendible bilinear forms. The case $n\ge 3$ is similar to the proof of Proposition~\ref{Para mas de 3 arruina}, using again that for bilinear forms extendibility is equivalent to domination.

Analogously, just as in the polynomial case, the results for $r$-integral and $r$-factorable multilinear forms (with the obvious definitions) can be deduced from the $r$-dominated case.

\bigskip

We end this note with some remarks on unconditionality for certain Banach operator ideals. We have seen that for two different Banach spaces, the lack of unconditionality in tensor products may fail. Therefore, it is reasonable to expect that, in order to obtain results of ``unconditionality destruction'' type, we must impose certain conditions to the involved spaces.

\smallskip

\textbf{The ideal of (p,q)-factorable operators \cite[Section 18]{DefFlo93}: }
Let $p,q \in [1, +\infty]$ such that $1/p + 1/q \geq 1$. An operator $T: E \to F$ is \textbf{$(p,q)$-factorable} if there are a finite measure $\mu$, operators $R \in \mathcal{L}\big( E,L_{q'}(\mu) \big)$ and $S \in \mathcal{L}\big( L_{p}(\mu), F''\big)$ such that $k_F \circ T = S \circ I \circ R$,
\begin{equation*}
\begin{array}{rcl}
E & \overset T \longrightarrow & F \overset{k_F} \hookrightarrow  F'' \\
_{R}\downarrow & & \; \; \; \; \nearrow_{S}\\
L_{q'}(\mu) & \underset I \longrightarrow & L_p(\mu),
\end{array}
\end{equation*}
where $I$ and $k_F$ are the natural inclusions.
We will denote the space of all such operators by $\Gamma_{p,q}(E,F)$.
For $T \in \Gamma_{p,q}(E,F)$, the (p,q)-factorable norm is given by $\gamma_{p,q}(T) = \inf \{ \|S\| \|I\| \|R\| \}$, where the infimum is taken over all such factorizations.

If $1/p + 1/q = 1$, $\Gamma_{p,q}$ coincides isometrically with the classical ideal $\Gamma_{p}$ of $p$-factorable operators \cite[Chapter 9]{DieJarTon95}.

\begin{example} \label{(p,q) factorables}
Let $E$ and $F$ be Banach spaces with unconditional basis such that $E'$ and $F$ have both finite cotype. Then, $\Gamma_{p,q}(E,F)$ does not have the Gordon-Lewis property.
\end{example}

\begin{proof}

By Theorem~\ref{tzafriri} we know that, for $r \in \{2, \infty\}$ and $s \in \{1,2\}$, $E$ and $F$ contain the uniformly complemented sequences $(\ell_r^m)_{m=1}^\infty$, $(\ell_s^m)_{m=1}^\infty$ respectively. This easily implies that $\Gamma_{p,q}(E,F)$ must contain the uniformly complemented sequence $(\Gamma_{p,q}(\ell_r^m,\ell_s^m))_{m=1}^\infty$.
Therefore, if show that $gl \big(\Gamma_{p,q}(\ell_r^m,\ell_s^m)\big) \to \infty$ as $m \to \infty$ we are done.

By \cite[17.10]{DefFlo93} we know that $(\Gamma_{p,q}, \gamma_{p,q})$ is a maximal operator ideal associated with the tensor norm $\alpha_{p,q}$ of Lapresté (see \cite[12.5]{DefFlo93} for definitions).
Thus,
$$\Gamma_{p,q}(\ell_r^m,\ell_s^m)=  \ell_{r'}^m \otimes_{\alpha_{p,q}} \ell_{s}^m.$$
Now by \cite[Exercise 31.2. (a)]{DefFlo93} we have $$gl\big( \Gamma_{p,q}(\ell_r^m,\ell_s^m) \big) = gl\big(\ell_{r'}^m \otimes_{\alpha_{p,q}} \ell_{s}^m \big) \asymp gl\big(\ell_{r'}^m \otimes_{/ \pi \s} \ell_{s}^m \big) $$ which goes to infinity as $m \to \infty$ (this is a direct consequence of the proof of Proposition~\ref{ideal destruction test} for $\mathcal{P}_e^2$ and Lemma~\ref{GL-chevet}).
\end{proof}

In particular, for $1 < r < \infty$ and $1 \leq s < \infty$ the spaces $\Gamma_{p,q}(\ell_r,\ell_s)$ and $\Gamma_{p,q}(c_0,\ell_s)$ do not have the Gordon-Lewis property. The case $r=\infty$ and $1 \leq s < \infty$ can be established just following the previous proof. In fact, proceeding as above and using \cite[Proposition 7]{Schutt}, something more can be stated:
For $2 \leq r \leq \infty$ and $1 \leq s \leq 2$, if $E$ and $F$ be Banach spaces such that $E$ contains the sequence $(\ell_r^m)_{m=1}^\infty$ uniformly complemented and $F$ contains the sequence $(\ell_s^m)_{m=1}^\infty$ uniformly complemented, then $\Gamma_{p,q}(E,F)$ does not have the Gordon-Lewis property. Note that in this case, we do not require that $E$ nor $F$ have unconditional bases.

\bigskip

\textbf{The ideal of (p,q)-dominated operators \cite[Section 19]{DefFlo93}:}
Let $p,q \in [1, +\infty]$ such that $1/p + 1/q \leq 1$. An operator $T: E \to F$ is \textbf{$(p,q)$ dominated} if for every $m \in \mathbb{N}$, $x_1, \dots, x_m \in E$ and $y_1', \dots, y_m' \in F'$ there exist a constant $C \geq 0$ such that:
$$ \ell_r (<y_k',Tx_k>) \leq C w_p(x_k) w_q(y_k'),$$
where $1/p + 1/q + 1/{r'} = 1$. We will denote the space of all such operators by $\mathcal{D}_{p,q}(E,F)$ with the norm $D_{p,q}(T)$ being the minimum of these $C$.

Equivalently, $T \in \mathcal{D}_{p,q}(E,F)$ if there are a constant $B \geq 0$ and probability measures $\mu$ and $\nu$ such that
$$ |<y',Tx>| \leq B \big( \int_{B_{E'}} |<x',x>|^p \mu(dx') \big)^{1/p} \big( \int_{B_{F''}} |<y'',y'>|^q \nu(dy'') \big)^{1/q},$$
holds for all $x \in E$ and $y' \in F'$, (replace the integral by $\| \; \|$ if the exponent is $\infty$).
In this case, the $(p,q)$-dominated norm of $T$, $D_{p,q}(T)$, is the infimum of the constants $B$ for which the previous inequality hold (see \cite[Corollary 19.2.]{DefFlo93}).

If $1/p + 1/q = 1$, $\mathcal{D}_{p,q}$ coincides isometrically with the classical ideal of $p$-dominated operators \cite[Chapter 9]{DieJarTon95}.

By \cite[Sections 17 and 19]{DefFlo93} we know that the ideal of $\mathcal{D}_{p,q}$ the adjoint of $\Gamma_{p',q'}$, the ideal of $(p',q')$-factorable operators. Using the duality that this implies on finite dimensional spaces, we can deduce:

\begin{example}
If $E$ and $F$ be Banach spaces with unconditional basis such that $E$ and $F'$ have both finite cotype. Then, $\mathcal{D}_{p,q}(E,F)$ does not have the Gordon-Lewis property.
\end{example}

As above, we can see that for $1 \leq r \leq 2$ and $2 \leq s \leq \infty$, if $E$ contains the sequence $(\ell_r^m)_{m=1}^\infty$ uniformly complemented and $F$ contains the sequence $(\ell_s^m)_{m=1}^\infty$ uniformly complemented, then $\mathcal{D}_{p,q}(E,F)$ does not have the Gordon-Lewis property.

We have, in particular,  that $1 \leq r < \infty$ and $1 \leq s \leq \infty$ the spaces $\mathcal{D}_{p,q}(\ell_r,\ell_s)$ and $\mathcal{D}_{p,q}(\ell_r,c_0)$ do not have the Gordon-Lewis property.

\bigskip
Let us give a procedure to obtain more examples: if $\mathcal{A}$ a Banach operator ideal and $\alpha$ is its associated tensor norm, by  $\mathcal{A}^{inj \; sur}$ we denote the maximal operator ideal associated to the norm $/ \alpha \s$ \cite[Sections 9.7 and 9.8]{DefFlo93}. Using the ideas of Example~\ref{(p,q) factorables} and the fact that $/ \alpha \s \leq / \pi \s$, we have:

\begin{example} \label{inj-sur}
Let $E$ and $F$ be Banach spaces with unconditional basis such that $E'$ and $F$ have both finite cotype. Then, $\mathcal{A}^{inj \; sur}(E,F)$ does not have the Gordon-Lewis property.
\end{example}

For example, let us consider $\mathcal A$ to be the ideal of $(p,q)$-factorable operators. An operator $T$ belongs to  $\Gamma_{p,q}^{inj \; sur}(E,F)$ if and only if there is a constant $C \geq 0$ such that for all natural numbers $m \in \mathbb{N}$, all matrices $(a_{k,l})$, all $x_1, \dots, x_m \in E$ and all $y_1', \dots, y_m' \in F'$

$$ \big| \sum_{k,l =1}^m a_{k,l} <y_k',Tx_l> \big| \leq C \|(a_{k,l}): \ell_{p'}^m \to \ell_{q}^m \| \ell_{p'}(x_l) \ell_{q'}(y_k').$$
In this case, $\gamma_{p,q}^{inj \: sur}(T) : = \min C$ (see \cite[Theorem 28.4]{DefFlo93}).

\end{document}